\documentclass[12pt]{amsart}

\usepackage{tikz}

\usepackage[utf8]{inputenc}

\usepackage{amsmath}
\usepackage{amssymb}
\usepackage{amsfonts}
\usepackage{amsthm}
\usepackage{graphicx}
\usepackage{enumerate}
\usepackage{a4wide}
\usepackage{tabularx}
\usepackage{subcaption}

\makeatletter
\@namedef{subjclassname@2010}{
  \textup{2010} Mathematics Subject Classification}
\makeatother

\usepackage{booktabs}

\usepackage{hyperref}

\newtheorem{theorem}{Theorem}

\newtheorem{corollary}[theorem]{Corollary}

\newtheorem{definition}[theorem]{Definition}
\newtheorem{question}[theorem]{Question}

\numberwithin{equation}{section}

\def\C{\mathbb{C}}

\def\Z{\mathbb{Z}}

\def\F{\mathbb{F}}
\def\al{\alpha}

\def\ll{\mathcal{L}}
\def\nn{\mathcal{N}}
\def\bb{\mathcal{B}}

\def\dd{\mathcal{D}}

\def\Rem{\mathcal{R}}
\def\G{\mathcal{G}}

\newcommand{\abs}[1]{\left|#1\right|}

\usepackage{color}

\begin{document}

\title[On certain multiples of Littlewood and Newman
polynomials]{On certain multiples of Littlewood and Newman
polynomials}
\author{P. Drungilas, J. Jankauskas, G. Junevi\v{c}ius, L. Klebonas, J. \v Siurys}
\address{Institute of Mathematics, Faculty of Mathematics and Informatics, Vilnius
University, Naugarduko 24, Vilnius LT-03225, Lithuania}
\email{pdrungilas@gmail.com}
\address{Mathematik und Statistik, Montanuniversit\"at Leoben, Franz Josef Stra\ss{}e 18, A-8700 Leoben, Austria}
\email{jonas.jankauskas@gmail.com}
\address{Institute of Mathematics, Faculty of Mathematics and Informatics, Vilnius
University, Naugarduko 24, Vilnius LT-03225, Lithuania}
\email{grintas@gmail.com}
\address{Institute of Mathematics, Faculty of Mathematics and Informatics, Vilnius
University, Naugarduko 24, Vilnius LT-03225, Lithuania}
\email{lukasklebonas@gmail.com}
\address{Institute of Mathematics, Faculty of Mathematics and Informatics, Vilnius
University, Naugarduko 24, Vilnius LT-03225, Lithuania}
\email{jonas.siurys@mif.vu.lt}

\thanks{This research of the first author was funded by a grant (No. S-MIP-17-66/LSS-110000-1274) 
from the Research Council of Lithuania. This research was performed 
in cooperation with the Department Mathematik und Statistik, Montanuniversit\"at Leoben}

\thanks{The research of J. Jankauskas was supported by FWF grant "Number Systems, Spectra and Rational Fractal Tiles" number M 2259-N35 under Lise Meitner Programme.}


\begin{abstract}
Polynomials with all the coefficients in $\{ 0,1\}$ and constant term 1 are called Newman polynomials, 
whereas polynomials with all the coefficients in $\{ -1,1\}$ are called Littlewood polynomials. By exploiting an algorithm  developed earlier, we determine 
the set of Littlewood polynomials of degree at most 12 which divide Newman 
polynomials. Moreover, we show that every Newman quadrinomial $X^a+X^b+X^c+1$, $15>a>b>c>0$, has a Littlewood multiple of smallest possible degree which can be as large as $32765$. 

%which decides whether a given monic integer polynomial with no roots on the unit circle $|z|=1$ has a non-zero multiple in $\Z[X]$ with coefficients in a finite set $\dd \subset \Z$,  for every 
%Borwein polynomial of degree at most 9 we determine whether it divides any Littlewood or Newman polynomial. In particular, we show that every Borwein polynomial of degree at most 8 which divides some Newman polynomial  divides some Littlewood polynomial as well. In addition to this, for every Newman polynomial of degree at most 11, we check whether it has a Littlewood multiple, extending the previous results of Borwein, Hare, Mossinghoff, Dubickas and Jankauskas.
\end{abstract}

\maketitle

\section{Introduction}

A polynomial
\begin{equation}\label{px}
P(X) = a_d X^d + a_{d-1}X^{d-1}+\dots+a_1X+a_0\in\Z[X]
\end{equation}
with a nonzero constant term $P(0)$ is called a \emph{Borwein polynomial}, if $a_j \in \{-1, 0, 1\}$ for each $0 \leq j \leq d$.  
To avoid trivialities,  we consider only polynomials with non-zero leading and constant terms $a_d \cdot a_0 \ne 0$. In such case, both $P(X)$ and its reciprocal polynomial $P^{*}(X) := X^dP(1/X)$ are of the same degree $d$. For example, $P(X) = X^5-X^3+X-1$  is a Borwein polynomial. 
%The set of all Borwein polynomials is denoted by $\bb$. 

Similarly, the polynomial $P(X)$ in \eqref{px} is called a \emph{Newman polynomial}, if all coefficients $a_j \in \{0, 1\}$ and $P(0)=1$. For instance, $P(X) = X^4+X+1$ is a Newman polynomial. The subset of $\Z[X]$ of all Newman polynomials is denoted by $\nn$. Finally, 
a Borwein polynomial $P(X)$ in \eqref{px} is called a \emph{Littlewood polynomial} if all 
of its coefficients are non-zero, i. e., $a_j \in \{-1, 1\}$ for each $0 \leq j \leq d$. 
For example,  $P(X) = X^4+X^3-X^2-X-1$  is a Littlewood polynomial. The set of all Littlewood polynomials is denoted by $\ll$. 
%One has trivial set relations $\nn \subset \bb$, $\ll \subset \bb$. 

We say that a polynomial $P(X)$ in \eqref{px} is a \emph{trinomial} if it has only three non-zero coefficients $a_j$, for $0 \leq j \leq d$. If the number of non-zero coefficients is four, $P(X)$ is called a \emph{quadrinomial}.

We say that a polynomial $P(X)$ has a Littlewood multiple if it divides some polynomial in the set $\ll$. We say that $P(X)$ has a Newman multiple if $P(X)$ divides some polynomial in $\nn$. When we need to restrict our attention only to polynomials of fixed degree, we use the subscript $d$ in $\nn_d$ and $\ll_d$ to denote the sets of Newman and Littlewood polynomials of degree $d$, respectively. We use the subscript ``$\leq d$'' to indicate the sets of polynomials of degree \emph{at most} $d$, that is
\[
\nn_{\leq d} = \bigcup_{j=0}^d \nn_j, \qquad \ll_{\leq d} = \bigcup_{j=0}^d \ll_j. %\qquad \bb_{\leq d} = \bigcup_{j=0}^d \bb_j.
\]
%Clearly, non-constant polynomials $P(X)$ with all non-negative coefficients 
% cannot have any positive real zeros $X \in [0, \infty)$. Newman 
%  polynomials are among such polynomials. To denote the subsets of 
%   Littlewood polynomials with no real positive zeros, we append 
%    the $"-``$ superscript, for instance,  $\ll^-$, $\ll_d^-$ and $\ll_{\leq d}^-$.  

Let $\mathcal{A} \subset \Z[X]$. We will employ the notation $\ll(\mathcal{A})$ to denote the set of 
polynomials $P(X) \in \mathcal{A}$ which divide some Littlewood polynomial. In the same way, denote by 
$\nn(\mathcal{A})$ the set of polynomials $P(X) \in \mathcal{A}$ which divide some Newman polynomial. 
In particular, the set $\nn_d\!\setminus\!\ll(\nn)$ consists of those Newman polynomials of degree $d$ 
that do not divide any Littlewood polynomial, whereas the set 
$\ll_{\leq d}\setminus \nn(\ll)$ consists of those  
Littlewood polynomials of degree at most $d$ that do not divide any Newman polynomial.

Let $\dd\subset\Z$ be a finite set. We call $\dd$ a \emph{digit set}. Consider 
the following question.

\begin{question}\label{p136}
Given a polynomial $P \in \Z[X]$, does there 
exist a nonzero polynomial with coefficients in $\dd$ which is divisible by $P$?
\end{question}

In \cite{DJS} we implemented an algorithm which answers the above question provided that the polynomial $P(X)$ has no roots on the unit circle $|z|=1$ in the complex plane. As far as we know the first instance of such an algorithm appears in the work of Frougny \cite{Fro} on the digit representations of numbers in algebraic integer bases produced by finite automata. Lau \cite{ksl} describes a slightly different version of this algorithm used to determine discreteness property of Bernoulli measures, focused to the case when $P(X)$ is a minimal polynomial of Pisot number. Subsequently, Hare and his coauthors \cite{BH,HM} used the same algorithm to verify the uniform discreteness property of the spectra of certain Pisot numbers. Stankov \cite{dst} considered spectra of non--Pisot algebraic integers. Akiyama, Thuswaldner and Za\"{i}mi \cite[Th.~3]{ATZ} treated Question \ref{p136} in the context of \emph{height reducing problem} and arrived to the algorithm that is essentially similar to the one that was devised by Frougny \cite{Fro}. Thus, Question~\ref{p136} has been fully answered for separable monic polynomials $P(X) \in \Z[X]$ with no roots on the unit circle. The version of this algorithm that was implemented in \cite{DJS} is also capable to deal with non-separable cases $P(X) \in \Z[X]$ (i.e., when $P(x) \in \Z[x]$ have multiple roots). 
It is not known whether the non-vanishing on $|z|=1$ can be dropped or not; it seems to be essential to the proof that the search terminates. In many practical cases this restriction can be circumvented (see \cite{DJS}, Sec.~4.1 and the note at the end of Sec.~3 therein).

In \cite{DJS}, we implemented this algorithm to answer 
Question~\ref{p136} for all Borwein polynomials of degree up to 9 and 
the digit sets $\dd = \{0,1\}$ and $\dd=\{-1,1\}$. For every 
Borwein polynomial of degree at most 9 it was determined whether 
it has a Littlewood multiple and whether it divides some Newman polynomial. Also, for every Newman polynomial 
$P(X)$ of degree at most $11$, it was checked whether $P(X)\in\ll(\nn)$. 
These computations extended the results 
previously obtained by Dubickas and Jankauskas \cite{DJ}, Borwein and Hare \cite{BH}, Hare and Mossinghoff \cite{HM}.

In the present paper we show that all Newman quadrinomials of degree at most 
15, with the small number of possible exceptions that are given in Table~\ref{NQ15}, have Littlewood multiples. Moreover, we prove that every polynomial $P\in\ll(\nn_{\leq 10})$ possesses 
a Littlewood multiple of smallest possible degree $\deg_2\widetilde{P}$ (see Section~\ref{lmld} for the details) as well as every quadrinomial $Q\in\nn_{\leq 14}$ does. We also compute the set $\nn(\ll_{\leq 12})$ and show that it contains exactly nine polynomials which are not products of cyclotomic polynomials.

This paper is organized as follows. Main results regarding  Littlewood multiples of Newman polynomials and those 
 regarding Newman multiples of Littlewood polynomials  are given in  Sections~\ref{lmnq} and \ref{nmlp}, respectively.
 In Section~\ref{lmld} we consider Littlewood multiples of least degree. A short description of the algorithm used in 
this paper is given in section~\ref{sda}. Finally, 
 our computations are described in Section~\ref{cmp}.

\section{Littlewood multiples of Newman quadrinomials}\label{lmnq}

In \cite{DJS, DJ} the sets $\ll(\nn_d)$,  $\nn_d\setminus \ll(\nn)$ have been completely determined for 
 $d\leq 11$. In particular, it was proved in \cite{DJ} that every Newman polynomial of degree at most 8 divides some Littlewood polynomial, that is, 
 $\ll(\nn_{\leq 8})=\nn_{\leq 8}$. On the other hand, $\#\nn_{9}\!\setminus\! \ll(\nn) = 18$, $\#\nn_{10}\!\setminus\! \ll(\nn) = 36$ and $\#\nn_{11}\!\setminus\! \ll(\nn) = 174$. It was also proved in \cite{DJ} that every Borwein trinomial 
\[
 X^{a} \pm X^{b} \pm 1, \quad 1 \leq b < a, \quad a, b \in\Z
\] 
(including Newman trinomials $X^a+X^b+1$)  has a Littlewood multiple, as well as certain Borwein quadrinomials 
$X^a\pm X^b\pm X^c\pm1$ do. Dubickas and Jankauskas \cite{DJ} asked: 
\begin{question}[\cite{DJ}]
Does there exist a Newman quadrinomial with no Littlewood multiple? 
Equivalently, does the set $\nn \!\setminus\! \ll(\nn)$ contain a quadrinomial?
\end{question}
In \cite{DJS}, this question was extended to include Borwein quadrinomials; 20 such quadrinomials  
in $\bb_{\leq 9}\!\setminus\!\ll(\bb)$ were found. 
However, none of them is a Newman quadrinomial (each of them has coefficients $-1$ and $1$). 

Our computations show that

\begin{theorem}\label{qt256}
Every quadrinomial $Q\in\nn_{\leq 15}$, 
except possibly for those that are given in Table~\ref{NQ15}, 
divides some Littlewood polynomial. 
Moreover, every
quadrinomial $Q\in\nn_{\leq 14}$ possesses a 
Littlewood multiple of smallest possible degree 
$\deg_2 \widetilde{Q}-1$ (this quantity will be defined in Section~\ref{lmld}).
\end{theorem}

\begin{figure}[h]
\captionsetup[subfigure]{labelformat=empty}
	\begin{subfigure}{.45\linewidth}
		\centering
      		\includegraphics[scale=0.35]{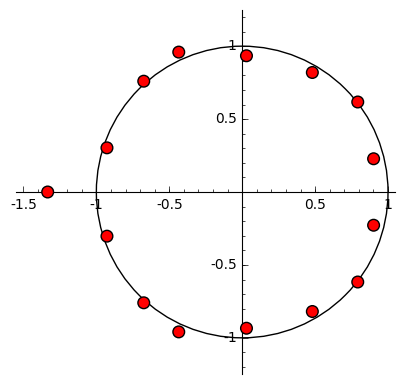}
      		\caption{$X^{15} + X^{14} + X^{10} + 1$.}
	\end{subfigure}
	\begin{subfigure}{.45\linewidth}
      		\centering
      		\includegraphics[scale=0.35]{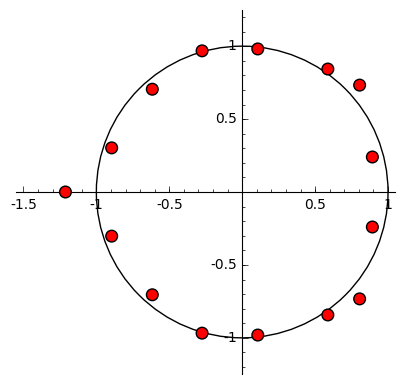}
      		\caption{$X^{15} + X^{12} + X^{10} + 1$.}
        \end{subfigure}
	\begin{subfigure}{.45\linewidth}
		\centering
		\includegraphics[scale=0.35]{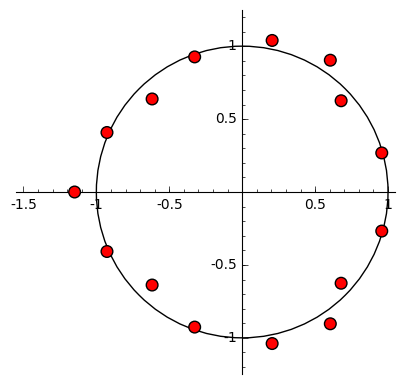}
		\caption{$X^{15} + X^{12} + X^{4} + 1$}
		\label{poly5}
        \end{subfigure}
        \begin{subfigure}{.45\linewidth}
        	\centering
        	\includegraphics[scale=0.35]{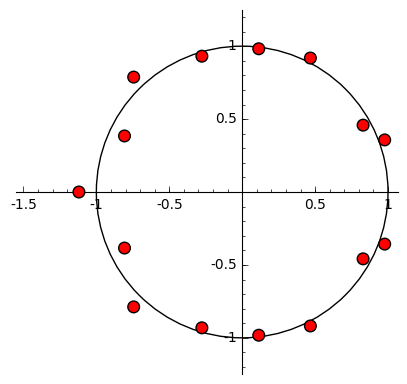}
        	\caption{$X^{15} + X^{8} + X^{6} + 1$}
	\end{subfigure}
	\caption{Complex roots of quadrinomials $Q(X)$ from Table~\ref{NQ15}.}\label{fig363}
\end{figure}

The second column in Table~\ref{NQ15} indicates the smallest possible degree 
of a Littlewood multiple (if it exists) of $Q\in\nn_{15}$. 
%It was proved in \cite{DJ} that if a  polynomial $P(X)\in\Z[X]$ divides a Littlewood 
%polynomial $L(X)$ then $\deg L\geq \deg_2 \widetilde{P}-1$, where the 
%number  $\deg_2 \widetilde{P}$ is defined in Section~\ref{lmld}.  
Figure~\ref{fig363} depicts the roots of polynomials from Table~\ref{NQ15}.

\begin{table}
\caption {Newman quadrinomials $Q\in\nn_{\leq 15}$ which are not known to have Littlewood multiples (reciprocals omitted).}\label{NQ15} 
\begin{tabular}{clc}
\toprule
&Quadrinomial $Q(X)$ & $\deg_2 \widetilde{Q}-1$\\
\midrule
&$X^{15} + X^{14} + X^{10} + 1$ & $10921$\\
%&$X^{15} + X^{5} + X + 1$ & $10921$\\
&$X^{15} + X^{12} + X^{10} + 1$ & $32765$\\
%&$X^{15} + X^{5} + X^{3} + 1$ & $32765$\\
&$X^{15} + X^{12} + X^{4} + 1$ & $31681$\\
%&$X^{15} + X^{11} + X^{3} + 1$ & $31681$\\
%&$X^{15} + X^{9} + X^{7} + 1$ & $32765$\\
&$X^{15} + X^{8} + X^{6} + 1$ & $32765$\\
\bottomrule
 \end{tabular}
 \end{table}

\section{Newman multiples of Littlewood polynomials}\label{nmlp}

We already mentioned, that the sets $\ll(\nn_d)$ have been completely determined for $d \leq 11$. Also, it was proved in \cite{DJS} (see Proposition~18) that if $P(X)$ is a product of cyclotomic polynomials and $P(1)\neq 0$ then $P(X)$ divides some Newman polynomial. Therefore Littlewood polynomials which are products of cyclotomic polynomials and do not vanish at $X=1$ have Newman multiples too.

Motivated by this result, we considered \emph{the complementary problem} of computing the sets $\nn(\ll_d)$ for small degrees $d$. By running the algorithm, described in Section~\ref{sda}, on every Littlewood polynomial $L(X) \in \ll_{\leq 12}$, we calculated explicitly all non-trivial elements of the set $\nn(\ll_{\leq 12})$.

\begin{theorem}
There are exactly $9$ monic Littlewood polynomials $P(X)$ of degree $\leq 12$ that possess Newman multiples and are not themselves the products 
of cyclotomic polynomials.  
%There are exactly nine monic Littlewood polynomials of degree $\leq 12$ which are not products 
%of cyclotomic polynomials and divide some Newman polynomial.  
\end{theorem}

\begin{figure}[h!]
\begin{subfigure}{0.30\textwidth}
\centering
\includegraphics[height=28mm]{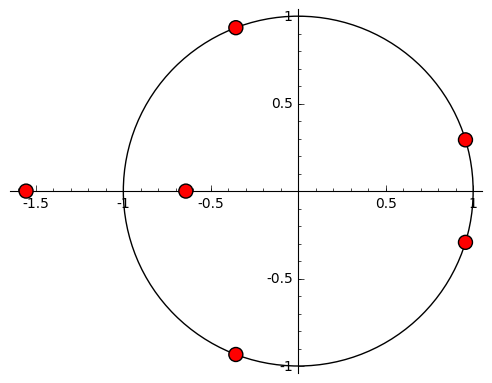} 
\caption*{\qquad$P_1(x)$}
\end{subfigure}\;\;
\begin{subfigure}{0.30\textwidth}
\centering
\includegraphics[height=28mm]{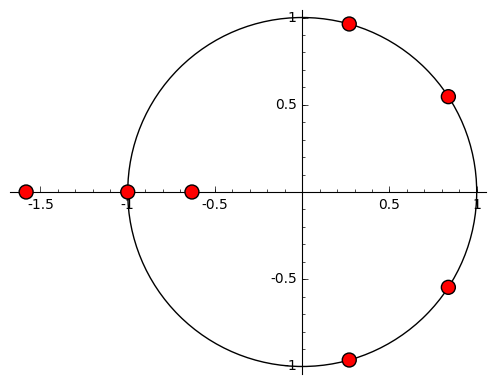} 
\caption*{\qquad$P_2(x)$}
\end{subfigure}\;\;
\begin{subfigure}{0.30\textwidth}
\centering
\includegraphics[height=28mm]{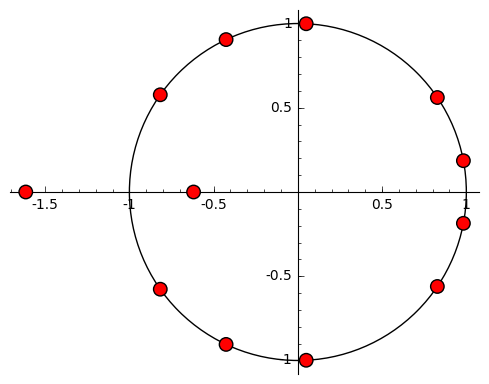}
\caption*{\qquad$P_7(x)$}
\end{subfigure}
%%%%%%%%%%%%%%%%%%
\begin{subfigure}{0.4\textwidth}
\centering 
 \includegraphics[height = 28mm]{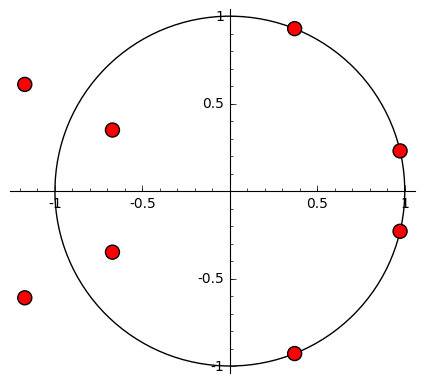}
\caption*{$P_3(x)$}
\end{subfigure}
\begin{subfigure}{0.4\textwidth}
\centering
\includegraphics[height=32mm]{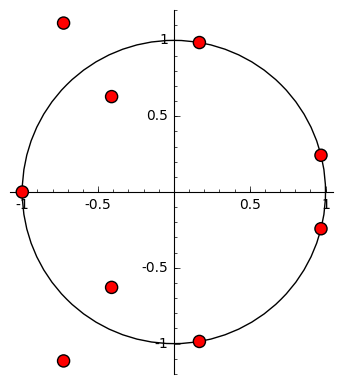}
\caption*{$P_4(x)$}
\end{subfigure}\\ %%%%%
\begin{subfigure}{0.4\textwidth}
\centering
\includegraphics[height=3.4cm]{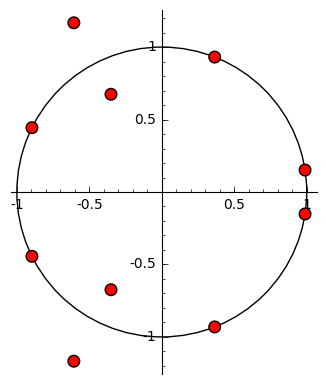} 
\caption*{$P_5(x)$}
\end{subfigure}
\begin{subfigure}{0.4\textwidth}
\centering
\includegraphics[height=3.4cm]{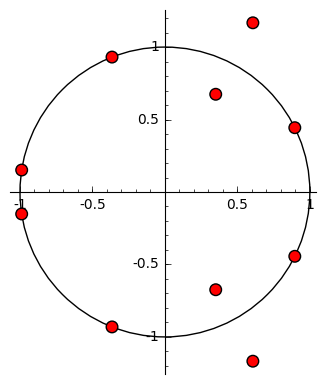}
\caption*{$P_6(x)$}
\end{subfigure}
%%%%%%%%%%%%%%%%%%
\begin{subfigure}{0.37\textwidth}
\centering
\includegraphics[height=34mm]{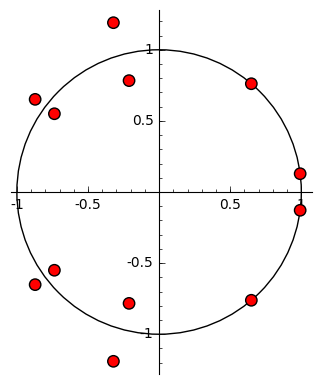} 
\caption*{$P_8(x)$}
\end{subfigure}\;\;
\begin{subfigure}{0.37\textwidth}
\centering
\includegraphics[height=34mm]{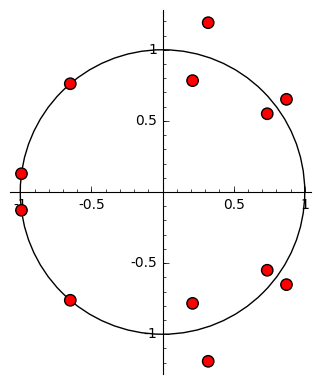}
\caption*{$P_9(x)$}
\end{subfigure}
\caption{Roots of polynomials $P_k(X)\in\nn(\ll_{\leq 12})$ from Table~\ref{nl12}}\label{lnp410}
\end{figure}

This result contrasts very sharply with the fact that \emph{most} Newman polynomials of small degree have Littlewood multiples.
Roots of these 9 polynomials are depicted in Figure~\ref{lnp410}.

Recall that  a real algebraic integer $\al > 1$ is called a \emph{Salem number} 
(see, e.g., \cite{Sa1, Sa2, Sa3}), if all other conjugates of $\al$ lie in the 
unit circle $|z| \leq 1$, with at least one conjugate on the unit circle $|z|=1$. 
Similarly, an algebraic integer $\al$ is called a \emph{negative Salem number} if $-\al$ is 
a Salem number. Finally, a \emph{complex Salem number} (see, e.g., \cite{Z}) is a nonreal algebraic 
integer $\al$ of modulus $\abs{\al}>1$ whose other conjugates, except for $\overline{\al}$, are of moduli $ \leq 1$, with at least one conjugate of
modulus $=1$. Note that 
the noncyclotomic parts (polynomials, obtained omitting their cyclotomic factors) of polynomials $P_1(X)$, $P_2(X)$ and $P_7(X)$ are 
minimal polynomials of negative Salem numbers, whereas the noncyclotomic 
parts of polynomials $P_3(X)$, $P_4(X)$, $P_5(X)$ and $P_6(X)$ are minimal polynomials 
of complex Salem numbers. As a result we have the following corollary.

\begin{table}[h]
\centering
\caption{Monic non-trivial elements of $\nn(\ll_{\leq 12})$}\label{nl12}
%\resizebox{\columnwidth}{!}{%
\begin{tabular}{ll}
\toprule
$k$ & $P_k(X)$\\
\midrule
$1$ & $X^6 + X^5 - X^4 - X^3 - X^2 + X + 1$ \\[0,1cm]
$2$ & $X^7 + X^6 - X^5 + X^4 + X^3 - X^2 + X + 1$ \\[0,1cm]
$3$ & $X^8 + X^7 - X^6 - X^5 + X^4 - X^3 - X^2 + X + 1$ \\[0,1cm]
$4$ & $X^9 + X^8 + X^7 - X^6 - X^5 - X^4 - X^3 + X^2 + X + 1$ \\[0,1cm]
$5$ & $X^{10} - X^9 + X^8 + X^7 - X^6 + X^5 - X^4 + X^3 + X^2 - X + 1$ \\[0,1cm]
$6$ & $X^{10} + X^9 + X^8 - X^7 - X^6 - X^5 - X^4 - X^3 + X^2 + X + 1$ \\[0,1cm]
$7$ & $X^{12} + X^{11} - X^{10} - X^9 - X^8 + X^7 + X^6 + X^5 - X^4 - X^3 - X^2 + X + 1$ \\[0,1cm]
$8$ & $X^{12} - X^{11} + X^{10} + X^9 - X^8 + X^7 + X^6 + X^5 - X^4 + X^3 + X^2 - X + 1$ \\[0,1cm]
$9$ & $X^{12} + X^{11} + X^{10} - X^9 - X^8 - X^7 + X^6 - X^5 - X^4 - X^3 + X^2 + X + 1$ \\[0,1cm]
\bottomrule
\end{tabular}
%}
\end{table}

%%%%%%%%%%%%%%%%%%

%%%%%%%%%%%%%%%%%%%%

\begin{corollary}
Every polynomial $P\in\nn(\ll_{\leq 12})$ that is not a product 
of cyclotomic polynomials, has a complex root on $|z|=1$ which is not a 
root of unity.  
\end{corollary}

\section{Littlewood multiples of minimal degree}\label{lmld}

Let $\F_2$ be a finite field with two elements. Consider a polynomial 
$p(X)\in\F_2[X]$ with a non-zero constant term $p(0)$ in 
$\F_2$. Then $p(X)$ has a 
unique representation as a product
\[
p(X) = (X+1)^m\prod_{j=1}^r\phi_j(X)^{m_j},
\]
where $m\geq 0$ and $\phi_j(X)\in\F_2[X]$ are irreducible polynomials of degree 
greater than or equal to 2 and multiplicity $m_j\geq 1$, $j=1,2,\dotsc, r$. 
The product is empty if $r=0$. Every polynomial $\phi_j(X)$ divides a unique cyclotomic polynomial $\Phi_{e_j}(X)$ of odd index $e_j$ (see, e.g., \cite[Section~4]{DJ}). Let $s$ be the least positive integer satisfying 
$2^s\geq \max\{ m+1,m_1,\dotsc,m_r\}$. Following \cite{DJ}, we define 
the number 
\[
\deg_2 p = 2^s\cdot \text{lcm}(e_1,\dotsc,e_r).
\]
Given a polynomial $P(X)\in\Z[X]$, denote by $\widetilde{P}(X)\in\F_2[X]$ its 
reduction modulo 2. In \cite{DJ}, it was shown that whenever $P(X)\in\Z[X]$ divides a Littlewood 
polynomial $L(X)$, then $\deg L = k\cdot\deg_2 \widetilde{P}-1$ for some $k\in\{1,2,\dotsc\}$. Therefore,  $\deg L\geq \deg_2 \widetilde{P}-1$. For example, 
$\deg_2 \widetilde{P}=32766$ for the polynomial $P(X)=X^{15} + X^{8}+X^6+1$ (see Table~\ref{NQ15}). 
Hence, if this particular quadrinomial has a Littlewood multiple $L(X)$, then 
 $\deg(L) \geq 32765$. Another example is the quadrinomial $Q(X) = X^{15} + X^{10}+X^8+1$: we have computed a Littlewood multiple of $Q(X)$ of degree $\deg_2\widetilde{Q}-1=32765$. However, in general there might be more than one monic Littlewood multiple of minimal degree. Indeed, we found four distinct monic Littlewood multiples of $P(X) = X^{12} + X^{11}+X^{10}+1$ of degree $\deg_2\widetilde{P}-1 = 1189$.   

Our computations yield the following result.

\begin{theorem}
Every polynomial $P(X)\in\ll(\nn_{\leq 10})$, except possibly for those given in 
 Table~\ref{neradau}, has a Littlewood multiple of smallest possible degree 
$\deg_2 \widetilde{P}-1$.
\end{theorem}

\begin{table}[h]
\centering
 \caption {Polynomials $P(X)\in\ll(\nn_{\leq 10})$ which are not known to have Littlewood multiples of degree $\deg_2 \widetilde{P}-1$ (reciprocals omitted).}\label{neradau}
\begin{tabular}{clc}
\toprule
$n$ & $P_n(X)$ & $\deg_2 \widetilde{P_n}$\\
\midrule
1&$X^{9} + X^{8} + X^{6} + X^{5} + X^{4} + X^{3} + X + 1$ & $60$\\
2&$X^{10} + X^{9} + X^{8} + X^{3} + X^{2} + 1$ & $1020$\\
%3&$X^{10} + X^{8} + X^{7} + X^{2} + X + 1$ & $1020$\\
3&$X^{10} + X^{9} + X^{7} + X^{6} + X^{5} + 1$ & $1022$\\
%5&$X^{10} + X^{5} + X^{4} + X^{3} + X + 1$ & $1022$\\
\bottomrule
 \end{tabular}
 \end{table}
 
\begin{figure}[h]
 \captionsetup[subfigure]{labelformat=empty}
	\begin{subfigure}{.45\linewidth}
		\centering
      		\includegraphics[scale=0.4]{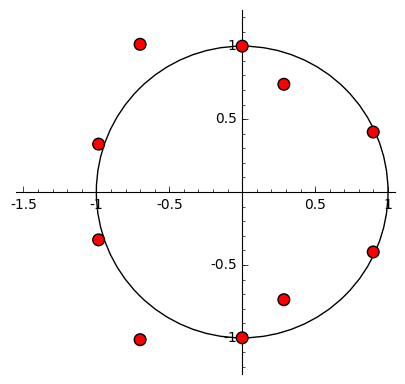}
      		\caption{$X^{10} + X^9 + X^8 + X^{3} + X^{2} + 1$}%\label{ex1P}
	\end{subfigure}
	\begin{subfigure}{.5\linewidth}
      		\centering
      		\includegraphics[scale=0.4]{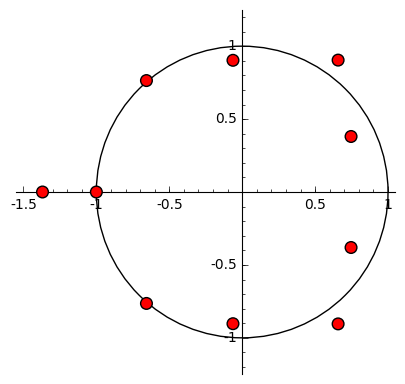}
      		\caption{$X^{10} + X^{9} + X^{7}+ X^{6} + X^{5} + 1$}%\label{ex1Q}
        \end{subfigure}
	\caption{Complex roots of $P_2(X)$ and $P_3(X)$ from Table~\ref{neradau}.}\label{ntp397}
\end{figure}

 The first polynomial in Table~\ref{neradau} is a product of cyclotomic polynomials:
 \[
% X^{9} + X^{8} + X^{6} + X^{5} + X^{4} + X^{3} + X + 1 = 
P_1(X)= (X + 1)^{3}  (X^{2} - X + 1)  (X^{4} - X^{3} + X^{2} - X + 1).
 \] 
Second and third polynomials in Table~\ref{neradau} factor in $\Z[X]$ as follows:
\begin{align*}
%X^{10} + X^{9} + X^{8} + X^{3} + X^{2} + 1 = 
P_2(X)&=(X^{2} + 1) (X^{8} + X^{7} - X^{5} + X^{3} + 1),\\
 %X^{10} + X^{9} + X^{7} + X^{6} + X^{5} + 1 =
P_3(X)&= (X + 1)  (X^{9} + X^{6} + X^{4} - X^{3} + X^{2} - X + 1). 
\end{align*}
Roots of $P_2(X)$ and $P_3(X)$ are depicted in Figure~\ref{ntp397}.

The corresponding result on Littlewood multiples of Newman quadrinomials was stated in Theorem~\ref{qt256} of Section~\ref{lmnq}.

\section{A short description of the algorithm}\label{sda}

%Let $\dd\subset\Z$ be a finite set, $B:=\max\{|d|: d\in\dd\}$.
We need the following definition.
\begin{definition}\label{def238}
Let $P(X)$ be a nonconstant polynomial with integer coefficients with no roots on the complex unit circle $|z|=1$. Suppose that the  
factorization of $P$ in $\C[X]$ is 
\[
P(X) = a\cdot (X - \al_1)^{e_1}(X - \al_2)^{e_2}\dotsb (X - \al_s)^{e_s},
\]
where $\al_1, \al_2, \dotsc, \al_s$ are distinct complex numbers and $e_j\geq 1$ for 
$j=1,2,\dotsc, s$. Let $B$ be arbitrary positive number. Define $\Rem(P, B)$ to 
be the set of all polynomials $R \in \Z[X]$, $\deg R < \deg P$, which, for each 
$j\in\{ 1,2,\dotsc, s \}$, satisfy the  
inequalities
%\begin{equation}
%\begin{eqnarray}\label{eq:237con}
%|R(\al_j)| \leq \frac{B}{||\al_j|-1|},\;\,|R'(\al_j)| \leq \frac{1!B}{||\al_j|-1|^2},\nonumber\\
%\dotsb\qquad\qquad\;\;\qquad\qquad\\
%|R^{(e_j-1)}(\al_j)| \leq \frac{(e_j-1)!B}{||\al_j|-1|^{e_j}}.\qquad\qquad\nonumber
%\end{eqnarray}
%\end{equation}
\begin{equation}\label{eq:237con}\mbox{\normalsize \(
|R(\al_j)| \leq \frac{B}{||\al_j|-1|},\;\,|R'(\al_j)| \leq \frac{1!B}{||\al_j|-1|^2},\dotsc, |R^{(e_j-1)}(\al_j)| \leq \frac{(e_j-1)!B}{||\al_j|-1|^{e_j}}.
\)}
\end{equation}
Here 
$R^{(k)}$ denotes the $k$th derivative of the polynomial $R$, and $R^{(0)}:=R$.
\end{definition}

Let $\dd\subset\Z$ be a finite set and let $\G=\G(P, \dd)$ be a directed graph whose vertices represent all the distinct polynomials 
$R \in \Rem(P, B)\cup \dd$, where $B=\max\{|b|: b\in \dd\}$.  
Two vertices corresponding to polynomials $R_i$ and $R_j$ are conected by a directed edge from $R_i$ to $R_j$, if $R_j \equiv X \cdot R_i + b \pmod{P}$ in $\Z[X]/(P)$ for some digit $b \in \dd$.

It was proved in \cite{ATZ, DJS, Fro} that, given a polynomial $P(X)\in\Z[X]$, which satisfies 
the conditions of Definition~\ref{def238}, the set $\Rem(P, B)$ is finite for every 
 $B>0$. Therefore, the graph $\G=\G(P, \dd)$ is finite for such polynomials $P$.

We restate Theorem~16 from \cite{DJS} without proof.

\begin{theorem}
Let $P \in \Z[X]$ be a monic polynomial with no roots on the complex unit circle $|z|=1$. Then $P$ divides an integer polynomial
\[
Q(X) = a_n X^n + a_{n-1} X^{n-1} + \dots  + a_1 X + a_0 \in \C[X]
\] with all the coefficients $a_j \in \dd$ and the leading coefficient $a_n \in \dd$, if and only if the graph $\G = \G(P, \dd)$ contains a path which starts at the remainder polynomial $R(X) = a_n$ and ends at $R(X) = 0$. The length of the path is $n$, where $n$ is the degree of $Q$.
\end{theorem}

The polynomial $Q$ with the coefficients in the set $\dd$ can be found by running any path finding algorithm on $\G=\G(P, \dd)$.

If $P(X)\in\Z[X]$ has unimodular roots that are not roots of unity, we exclude them when checking the inequalities \eqref{eq:237con} of Definition~\ref{def238} and continue building the graph $\G(P, \dd)$. Even in cases where $\G(P, \dd)$ is infinite, the component that is accessible from $R_0 = a_n$ might be finite, or this component might be simple enough so an existing path from $R_0$ to $R=0$ can be found within the reasonable time. Thus, in practice,  it is often the case that Question~\ref{p136} can be answered for such polynomials.

\section{Computations}\label{cmp}

We implemented the algorithm described in Section~\ref{sda} in C using Arb \cite{arb} library for arbitrary-precision floating-point 
ball arithmetic and ran it on the SGI Altix 4700 server at Vilnius University.  OpenMP \cite{omp} was used for the multiprocessing. 

For every $P \in \ll_{\leq 12}$, we calculated whether it divides some Newman polynomial. In the same way, for every Newman quadrinomial $Q\in\nn_{\leq 15}$, except for those given Table~\ref{NQ15}, we found a Littlewood multiple of $Q$. A slightly modified code was used to search for the smallest degree Littlewood multiples of Newman polynomials 
$P\in\ll(\nn_{10})$ and Newman quadrinomials $Q\in\nn_{\leq 15}$ (see Section~\ref{lmld}).  We will briefly explain 
how these calculations were organized. 

By Proposition~18, Section~4.1 of \cite{DJS}, one can omit cyclotomic factors of  
$P(X)\in\Z[X]$ when searching for its Littlewood multiples: if $\Phi_n(X)$ is the 
$n$-th cyclotomic polynomial, then 
$\Phi_n P \in\ll(\Z[X])$ if and only if $P\in\ll(\Z[X])$. In particular, every polynomial  that 
is a product of cyclotomic polynomials has a Littlewood multiple. For every quadrinomial $Q\in\nn_{\leq 15}$ we omited all of its cyclotomic factors and ran 
the algorithm to check whether $Q$ has a Littlewood multiple. Similarly, 
we ran the algorithm for every $P\in\ll(\nn_{\leq 10})$. Only for polynomials given in Table~\ref{NQ15} 
and Table~\ref{neradau} the computations did not terminate normally due to limits on our server machine.

Similarly, one can omit cyclotomic factors, except for $\Phi_1(X)=X-1$,  when searching for Newman multiples: for $n>1$, $\Phi_n P \in\nn(\Z[X])$ if and only if 
$P \in\nn(\Z[X])$. Also, real roots of  Newman polynomials must be negative, therefore $P \in\nn(\Z[X])$ implies 
that $P$ has no positive real roots. We used SAGE \cite{SAGE} to filter out all such polynomials from $\ll_{\leq 12}$ and then ran the algorithm on the  noncyclotomic parts of the remaining polynomials to check whether they have Newman multiples.

Some of the polynomials in sets $\ll_{\leq 12}$ and $\nn_{\leq 10}$ have unimodular roots that are not roots of unity. However, we succeeded in determining whether these polynomials belong to $\ll(\bb)$ and $\nn(\bb)$ (see Section \ref{sda} above). 

Note that in the inequalities \eqref{eq:237con} of Section \ref{sda}, one has the upper bound $B=\max\{|a|\mid a \in \dd\}=1$ in case of Littlewood and Newman multiples. To fasten the search, a new variable $\delta$, $0\leq \delta <1$ was introduced, replacing the bound $B$ by $B - \delta$. This eliminates some of the vertices in the original graph $\G(P, \dd)$. We start  
with the initial value $\delta = 0.95$.  If a Littlewood (or Newman) multiple is found, then we are done.  Otherwise we decrease $\delta$ by $0.05$ and try again. For polynomials in $\nn\!\setminus\!\ll(\nn)$ and $\ll\!\setminus\!\nn(\ll)$ the variable $\delta$ always reaches the value $\delta=0$ in order to construct the 
full graph $\G(P, \dd)$.
 
%\begin{figure}[h]
%    \caption{Distribution of noncyclotomic factors $F(X)$ of polynomials from $\bb_{\leq 9}$
%    such that $F(X)\in\ll(\Z[X])$.}\label{fdpp}
%      \centering
%      \includegraphics[scale=0.33]{littlewood}
%\end{figure}

All computations took approximately $49$ hours of CPU time. 
The maximum recursion depth reached by the algorithm was $98\,936$. It took about $90$ seconds of CPU time to run our 
algorithm to find the Littlewood multiple $L(X)$ of   
$P(X)=X^{15} + X^{10} + X^{8} + 1$ of smallest possible degree 
$\deg L = \deg_2\widetilde{P}-1=32\,765$; $283\,501$ vertices of the 
graph $\G(P,\{-1,1\})$ were constructed.  
On the other hand, it took $5960$ seconds of CPU time to find a Littlewood multiple $L(X)$ of  $P(X)=X^{14} - X^{12} + X^{3} +1$ of smallest possible degree 
$\deg L = \deg_2\widetilde{P}-1=16\,381$; more than one million vertices of the 
graph $\G(P,\{-1,1\})$ were constructed. For every polynomial in Table~\ref{NQ15} and Table~\ref{neradau} it took 2 hours of CPU time before the computation was interrupted due to limitations of computational resources.

\end{document}